\documentclass[journal]{IEEEtran}

\renewcommand{\baselinestretch}{1.0}
\renewcommand{\theenumii}{\theenumi.\arabic{enumii}}

\usepackage{amssymb,amsmath,graphicx,color,psfrag}
\usepackage{enumerate}
\usepackage{url,adjustbox}
\usepackage{graphicx}
\usepackage[utf8]{inputenc}
\usepackage{color,soul}
\usepackage[T1]{fontenc}
\usepackage{upgreek}
\usepackage{multicol}
\usepackage{multirow}
\usepackage{textcomp}
\usepackage{graphicx}
\usepackage{bm}
\usepackage{epstopdf}
\usepackage{amssymb}
\usepackage{amsmath}
\usepackage{psfrag}
\usepackage{amsthm}
\usepackage[tight,footnotesize]{subfigure}
\usepackage{cite}
\usepackage{enumerate}
 
\usepackage{mathtools}
\usepackage{array}
\usepackage{mathabx}
\usepackage{booktabs}
\usepackage{amsmath}
\usepackage{mathrsfs}
\usepackage{array}
\usepackage[tight,footnotesize]{subfigure}
\usepackage{algorithmic,algorithm}
\usepackage{multirow}
\usepackage{bm}
\usepackage{mathtools}
\usepackage{soul}
\usepackage{url}

\usepackage{upgreek}
\usepackage{epstopdf}
\usepackage{mathrsfs}
\usepackage{array}
\usepackage{multicol}
\usepackage[tight,footnotesize]{subfigure}
\usepackage{multirow}
\usepackage{bm}
\usepackage{soul}
\usepackage{savesym}

\newcounter{parcount}
\newcounter{parcountA}

\theoremstyle{plain}
\newtheorem{assumption}{Assumption}
\newtheorem{theorem}{Theorem}
\newtheorem{lemma}{Lemma}
\theoremstyle{definition}

\newtheorem{example}{Example}
\newtheorem{remark}{Remark}

\makeatletter
\newcommand{\vast}{\bBigg@{3.2}}
\newcommand{\Vast}{\bBigg@{5}}

\makeatother

\graphicspath{{figures/}}

\usepackage{etoolbox}
\makeatletter 
\pretocmd\@bibitem{\color{black}\csname keycolor#1\endcsname}{}{\fail}
\newcommand\citecolor[1]{\@namedef{keycolor#1}{\color{blue}}}
\makeatother
 
\begin{document}

\title{On Symmetric Gauss-Seidel ADMM Algorithm \\for $\mathcal{H}_\infty$ Guaranteed Cost Control with\\ Convex Parameterization}

\author{Jun Ma, 
	Zilong Cheng, Xiaoxue Zhang, 
	Masayoshi Tomizuka, \IEEEmembership{Life Fellow,~IEEE,}
	and
	Tong Heng Lee 
	\thanks{J. Ma is with the Robotics and Autonomous Systems Thrust, The Hong Kong University of Science and Technology (Guangzhou), Guangzhou, China, with the Department of Electronic and Computer Engineering, The Hong Kong University of Science and Technology, Hong Kong SAR, China, and also with the HKUST Shenzhen-Hong Kong Collaborative Innovation Research Institute, Futian, Shenzhen, China (e-mail: jun.ma@ust.hk).}
\thanks{Z. Cheng, X. Zhang, and T. H. Lee are with the Department of Electrical and Computer Engineering, National University of Singapore, Singapore 117583 (e-mail: zilongcheng@u.nus.edu; xiaoxuezhang@u.nus.edu; eleleeth@nus.edu.sg).}
\thanks{M. Tomizuka is with the Department of Mechanical Engineering, University of California, Berkeley, CA 94720 USA (e-mail: tomizuka@berkeley.edu).}
\thanks{This work has been submitted to the IEEE for possible publication. Copyright may be transferred without notice,	after which this version may no longer be accessible.}
}

\maketitle

\begin{abstract}
	This paper involves 
	the innovative
	development of a symmetric Gauss-Seidel ADMM algorithm 
	to solve the $\mathcal{H}_\infty$ guaranteed cost control problem. 
	In the presence of parametric uncertainties, 
	the $\mathcal{H}_\infty$ guaranteed cost control problem 
	generally leads to the large-scale optimization. This is due to the exponential growth of the number of the extreme systems involved with respect to the number of parametric uncertainties. In this work, through a variant of the Youla-Kucera parameterization, 
	the stabilizing controllers are parameterized in a convex set; 
	yielding the outcome 
	that the $\mathcal{H}_\infty$ guaranteed cost control problem 
	is converted to a convex optimization problem. 
	Based on an appropriate re-formulation using the Schur complement, 
	it then renders possible
	the use of the ADMM algorithm with symmetric Gauss-Seidel backward and forward sweeps. 
	Significantly, this approach alleviates the often-times prohibitively heavy computational burden 
	typical in 
	many $\mathcal H_\infty$ optimization problems 
	while exhibiting good convergence guarantees, 
	which is particularly essential for the related large-scale optimization procedures involved. 
	With this approach, the desired robust stability is ensured, 
	and the disturbance attenuation is maintained at the minimum level 
	in the presence of parametric uncertainties. 
	Rather importantly too, 
	with the attained effectiveness,
	the methodology thus evidently possesses extensive applicability 
	in various important controller synthesis problems, 
	such as decentralized control, sparse control, and output feedback control problems.
\end{abstract}


\section{Introduction}
Robust control theory typically
investigates the effect of disturbances, noises, 
and uncertainties on system performance; 
and continued great efforts have been devoted 
to robust stabilization and robust performance 
in the literature~\cite{kai2009robust,wang2016approximate,zhu2018output,wang2020global,ma2020small,chen2020novel,chen2020optimal}. 
Quite remarkably, 
several significant results~\cite{petersen1986riccati,petersen1987stabilization,bernussou1989linear} 
have been reported 
which
relate the notion of quadratic stabilization to robust stabilization 
for a class of uncertain linear systems, 
and by this concept, the stability of an uncertain system is established 
with a quadratic Lyapunov function. 
On the other hand, 
it is also the case that
$\mathcal H_\infty$ control 
is commonly and extensively used 
to attenuate the effect of disturbances 
on the system performance~\cite{doyle1989state,xie1992robust}. 
Additionally,
it is further known and shown in~\cite{zhou1988algebraic}
that a certain type of quadratic stabilization problem 
can be essentially expressed as an $\mathcal{H}_\infty$ control problem, 
where a Riccati inequality condition relates the determination of 
a stabilizing feedback gain 
that imposes a suitable $\gamma$ disturbance attenuation level~\cite{chang2020robust,wu2020fuzzy}. 
Also notably, 
the problem of finding the minimal disturbance attenuation level 
is recognized as an important and commonly-encountered problem,
and stated
as the optimal $\mathcal{H}_\infty$ control problem. 
On this, it is 
additionally noteworthy that the work
in~\cite{scherer1989h}
shows
that the problem can be tackled by an iterative algorithm 
based on the Riccati inequality condition. 
However here, nonlinear characteristics of the Riccati inequality condition 
typically result in 
significant complexity and difficulty 
in
obtaining the optimal gain and disturbance attenuation level.
Moreover, $\mathcal{H}_\infty$ filtering has been widely studied for state estimation~\cite{qi2019event,qi2021h}. It is remarkable that $\mathcal{H}_\infty$ filtering allows the system disturbances to be unknown, and uncertainties are tolerated in the system.

As considerable efforts have been made on the well-known Youla-Kucera parameterization 
(also known as $Q$-parameterization) for the determination 
of stabilizing  controllers~\cite{qi2004structured,furieri2019input,lin2019convex}, 
one may thus think about borrowing this idea to solve 
the $\mathcal{H}_\infty$ optimal control problem in the presence of parametric uncertainties. 
However, the derivations of the classical Youla-Kucera parameterization results 
rely on the fact that the plant is linear with no parametric uncertainty, 
and the order of the controller depends on the order of the plant model and that of $Q$. 
Alternative parameterization techniques based on the positive real lemma 
and the bounded real lemma~\cite{ma2019onrobust,ma2020robust} 
have also been proposed to deal with parametric uncertainties. 
However, as the required transfer function representation there
results 
in reduced stability in numerical computations, 
and high computational cost also incurs; 
it is not considered as a preferable choice for many practical applications. 
Hence, several other parameterization methods 
are presented instead in a state-space framework, 
for example~\cite{geromel1991convex}. 
In essence, these techniques are considered as 
variants of the Youla-Kucera parameterization, 
but with more flexibility to deal with the structural constraints and parametric uncertainties. 
Regarding the nonlinear constraints existing 
in such a parameterization (also noted to be convex in~\cite{geromel1991convex}), 
outer linearization is necessary for polyhedral approximation 
during iterative refinement~\cite{bertsekas2011unifying}. Similarly, the cutting-plane method is presented in~\cite{puttannaiah2015generalized,puttannaiah2016generalized} to solve generalized $\mathcal{H}_\infty$ control problems.
This technique could be effective in some scenarios, 
but there could be other certain scenarios such as high-dimension systems 
or uncertain systems with a large volume of parametric uncertainties.
For example, in~\cite{bergeling2020h}, a large-scale $\mathcal{H}_\infty$ optimal control problem is simplified in
a way that enables sparse solutions and efficient computation. In~\cite{zheng2017platooning}, a distributed $\mathcal{H}_{\infty}$ controller is developed for large-scale platoons.
In the scenarios that high-dimension optimization is involved (which resulted from rather numerous extreme systems involved computationally), 
the convergence rate of the cutting-plane method can be unacceptably slow; 
and in some cases, the optimization process could even terminate abruptly with unsuccessful outcomes. 
This is because a sizeable number of cutting planes needs to be added computationally at each iteration,
and in difficult scenarios, the optimization process can thus become unwieldy. 
It is also worth mentioning that this method can only guarantee 
the so-called $\epsilon$-optimality 
because the constraints are typically not exactly satisfied 
but violated by certain small values. 
Therefore, such a situation causes deviations from the ``true" optimal result, 
and consequently
the desired robustness is not perfectly guaranteed,
and particularly so if the parametric uncertainties are significant.

Because of the typical computational burden arising from 
the growth of system dimensions and parametric uncertainties, 
several advanced optimization techniques are presented in the more recent works. 
Hence for large-scale and nonlinear optimization problems, 
the
alternating direction method of multipliers (ADMM)~\cite{boyd2011distributed,sun2015convergent,makhdoumi2017convergence,du2019admm} 
has attracted considerable attention from researchers, 
and is widely used in various areas 
such as statistical learning~\cite{boyd2011distributed}, distributed computation~\cite{chen2017distributed,wu2018distributed,wei2019distributed}, and multi-agent systems~\cite{peng2018consensus,yang2018output}.
ADMM demonstrates high efficiency in the determination 
of the optimal solution to many challenging problems.
Remarkably too, some of these challenging optimization problems 
cannot even be solved by the existing conventional gradient-based approaches, 
and in these, ADMM demonstrates its superiority. 
Nevertheless, the conventional ADMM methodology only ensures 
appropriate convergence 
with utilization of 
a two-block optimization structure, 
and this constraint
renders a serious impediment to practical execution~\cite{chen2016direct}. 
To cater to this deficiency, the symmetric Gauss-Seidel technique 
can be used to conduct the ADMM optimization serially~\cite{adegbege2016gauss,chen2017efficient}, 
which significantly improves the feasibility of the ADMM in many large-scale optimization problems. 
Although these methodologies are reasonably well-established, 
nevertheless only rather generic procedures are given at the present stage. 
Therefore, it leaves an open problem on 
how to apply these advanced optimization techniques in control problems 
such that these methodologies can be extended beyond the theoretical level.

It is also rather essential at this point to note that in the presence of significant parametric uncertainties, 
the $\mathcal H_\infty$ optimization problem 
is usually of the large-scale type
(because of the exponential growth of the number of extreme systems involved computationally
with respect to the number of parametric uncertainties; 
and each of these extreme systems has a one-to-one correspondence 
to an inequality constraint to ensure the closed-loop stability). 
Therefore, our contributions are summarized as follows. Here, we propose a novel optimization technique 
to solve the resulting large-scale $\mathcal{H}_\infty$ guaranteed cost control problem resulting from parametric uncertainties, 
where the stabilizing controllers 
are characterized 
by an appropriate convex parameterization 
(which will be described and established analytically).
Firstly, we construct a convex set 
such that all the stabilizing controller gains 
are mapped onto the parameter space, 
and the desired robust stability is then attained 
with the optimal disturbance attenuation level 
in the presence of convex-bounded parametric uncertainties. With this parameterization technique, the parametric uncertainties can be suitably considered in the problem formulation.
Secondly, a suitably interesting problem re-formulation based on the Schur complement 
facilitates 
the use of the symmetric Gauss-Seidel ADMM algorithm. Comparing with the methods in the existing literature (as described previously), this approach alleviates the ofter-times prohibitively heavy computational burden typically in many large-scale $\mathcal H_\infty$ control problems.

The remainder of this paper is organized as follows. 
In Section II, 
the 
optimal $\mathcal{H}_\infty$ controller synthesis with convex parameterization is provided. 
Section III presents the symmetric Gauss-Seidel ADMM algorithm to solve the $\mathcal{H}_\infty$ guaranteed cost control problem. 
Then, to validate the proposed algorithm, 
appropriate illustrative examples are given in Section IV with simulation results. 
Finally, pertinent conclusions are drawn in Section V.

\section{Optimal $\mathcal{H}_\infty$ Controller Synthesis by Convex Parameterization}
\textit{Notations:} $\mathbb R^{m\times n}$ ($\mathbb R^{n}$) denotes the real matrix with $m$ rows and $n$ columns ($n$ dimensional real column vector). $\mathbb S^{n}$ ($\mathbb S^{n}_{+}$) denotes the $n$ dimensional (positive semi-definite) real symmetric matrix, and $\mathbb S^{n}_{++}$ denotes the $n$ dimensional positive definite real symmetric matrix. The symbol $A \succ 0$ ($A \succeq 0$) means that the matrix $A$ is positive definite (positive semi-definite). $A^T$ ($x^T$) denotes the transpose of the matrix $A$ (vector $x$). $I_n$ ($I$) represents the identity matrix with a dimension of $n\times n$ (appropriate dimensions). The operator $\operatorname{Tr}(A)$ refers to the trace of the square matrix $A$. The operator $\langle A, B \rangle$ denotes the Frobenius inner product i.e. $\langle A,B\rangle= \operatorname{Tr}\left(A^TB\right)$ for all $A,B \in \mathbb R^{m\times n}$. The norm operator based on the inner product operator is defined by $\|x\|=\sqrt{\langle x,x\rangle}$ for all $x\in \mathbb R^{m\times n}$. $\| H(s)\|_\infty$ represents the $\mathcal{H}_\infty$-norm of $H(s)$. The operator $\operatorname{vec}(\cdot)$ denotes the vectorization operator that expands a matrix by columns into a column vector. The symbol $\otimes$ denotes the Kronecker product. $\sigma_\textup{M}(\cdot)$ returns the maximum singular value.

Consider a linear time-invariant (LTI) system
\begin{IEEEeqnarray*}{rCl}
\dot x &=& A x + B_2 u + B_1 w \IEEEyesnumber\IEEEyessubnumber\label{eq:ssm1} \\
z &=& C x+ D u   \IEEEyessubnumber\label{eq:ssm2}\\
u &=& -K x, \IEEEyessubnumber\label{eq:ssm3}
\end{IEEEeqnarray*}
with $x(0)=x_0$, $x\in \mathbb{R}^{n}$ is the state vector, $u\in \mathbb{R}^{m}$ is the control input vector, $w\in \mathbb{R}^{l}$ is the exogenous disturbance input, $z\in \mathbb{R}^{q}$ is the controlled output vector, $K \in \mathbb{R}^{m \times n}$ is the feedback gain matrix.
As a usual practice, Assumption~\ref{assu:usual} is made.
\begin{assumption}~\label{assu:usual}
	$[A, B_2]$ is stabilizable with disturbance attenuation $\gamma$, $[A, C]$ is observable, $C^T D=0$, and    $D^T D\succ 0$.
\end{assumption}

Denote $A_{c} = A- B_2 K$ and $C_{c} = C- D K$, the transfer function from $w$ to $z$ is given by
\begin{IEEEeqnarray}{rl}
H(s)=C_{c} (sI_n-A_{c})^{-1} B_1,  \label{eq:transfer function}
\end{IEEEeqnarray}
and the $\mathcal{H}_\infty$-norm is defined as \begin{IEEEeqnarray}{rl}\|H(s)\|_\infty=\sup \limits_\omega \,\sigma_\textup{M}[H(j\omega)].\end{IEEEeqnarray} 



It is worth stating that the objective of the optimal $\mathcal{H}_\infty$ control problem is to design a feedback controller that minimizes the $\mathcal{H}_\infty$-norm while maintaining the closed-loop stability. When the plant is affected by parametric uncertainties, the minimization of the upper bound to the $\mathcal{H}_\infty$-norm under all feasible models is known as the $\mathcal{H}_\infty$ guaranteed cost control problem. Note that in this work, $\gamma$-attenuation means that the $\mathcal{H}_\infty$-norm of $H(s)$ is bounded by $\gamma$, i.e., $\|H(s)\|_\infty\leq \gamma$.

In this work, for brevity, we define $p=m+n$ and $r=m+2n$. Then, the following extended matrices are introduced to represent the open-loop model \eqref{eq:ssm1}-\eqref{eq:ssm2}:
\begin{gather}
F=\begin{bmatrix}
A & -B_2\\
0 & 0
\end{bmatrix} \in \mathbb{R}^{ p \times  p}, \,
G=\begin{bmatrix}
0\\I_m
\end{bmatrix} \in \mathbb{R}^{ p \times m}, \nonumber\\
Q=\begin{bmatrix}
B_1 B_1^T &0\\0&0
\end{bmatrix} \in \mathbb{S}^{ p \ } , \, R=\begin{bmatrix}
C^T C &0\\0& D^T D
\end{bmatrix} \in \mathbb{S}^{ p  }.
\end{gather}
Also, define the matrix
\begin{IEEEeqnarray}{rl}
W=W^T=\begin{bmatrix}
W_1 & W_2 \\
W_2^T & W_3
\end{bmatrix},  \label{eqn:partition} \end{IEEEeqnarray} where $W_1 \in \mathbb{S}^{ n  }_{++}$, $W_2 \in \mathbb{R}^{ n \times m}$, $W_3 \in \mathbb{S}^{m  }$, and then define the matrical function \begin{IEEEeqnarray}{rl}\Theta(W, \mu)= FW+WF^T+WRW+\mu Q, \end{IEEEeqnarray} with $\mu=1/{\gamma}^2$~\cite{peres1994optimal}. Similarly, $\Theta(W, \mu)$ is partitioned as \begin{IEEEeqnarray}{rl}~\label{eq:theta_mu}
\Theta(W, \mu)= \begin{bmatrix} \Theta_{1}(W, \mu) & \Theta_{2}(W) \\ \Theta_{2}^T(W) & \Theta_{3}(W) \end{bmatrix},
\end{IEEEeqnarray}
with $\Theta_{1}(W, \mu) \in \mathbb{S}^{ n }, \Theta_{2}(W) \in \mathbb{R}^{ n \times m}, \Theta_{3}(W) \in \mathbb{S}^{  m}$.

The following theoretical developments generalize the results in~\cite{peres1994optimal,ma2021optimal,ma2021convex1}.
\begin{theorem} \label{thm:convex}
Define the set $\mathscr{C} = \{(W, \mu): W=W^T \succeq 0, \Theta_1(W, \mu) \preceq 0, \mu > 0$\}. Then the following statements hold:
\begin{enumerate}[(a)]	
	\item $\mathscr{C}$ is a convex set.
	\item Any $(W, \mu) \in \mathscr{C}$ generates a stabilizing gain $K=W_2^T W_1^{-1}$ that guarantees $\|H(s)\|_\infty \leq \gamma$ with $\gamma=1/\sqrt \mu>0$.
	\item 	At optimality, $(W^*,\mu^*)=\textup{argmax}\{\mu:(W, \mu)\in \mathscr{C}\}$ gives the optimal solution to the optimal $\mathcal{H}_\infty$ control problem, with $K^*={W_2^*}^{T} {W_1^*}^{-1}$ and $\|H(s)\|_\infty^*=\gamma^*=1/\sqrt{\mu^*}$.
\end{enumerate}
\end{theorem}

\noindent{\textbf{Proof of Theorem \ref{thm:convex}:}}
For Statement (a), the convexity of $\mathscr{C}$ can be proved as follows: first, the set of all positive semi-definite $W$ is a convex cone; second, for $\Theta(W)$: because $FW+WF^T$ is affine with $W$ and $\mu Q$ is linear with $\mu$; then, it remains to prove that $WRW$ is convex. Notably here, we only need to prove the convexity instead of the strong convexity. Take symmetric positive semi-definite matrices $W^1$ and $W^2$, then we have $\alpha W^1+(1-\alpha) W^2$ is symmetric, with $\alpha \in [0,1]$. Assume $\alpha W^1+(1-\alpha) W^2 \succeq 0$, we have
\begin{IEEEeqnarray}{rCl}
&&WRW  \nonumber\\
&=&\left[\alpha W^1+(1-\alpha) W^2\right] R \left[\alpha W^1+(1-\alpha) W^2\right] \nonumber\\
&=& \alpha^2 W^1 R W^1 + (1-\alpha)^2 W^2 R W^2 +2\alpha (1-\alpha)W^1 R W^2 \nonumber\\
&=& \alpha W^1RW^1+(1-\alpha)W^2RW^2\nonumber\\&& +\alpha (\alpha-1) (W^1 R W^1 +W^2 R W^2  - 2W^1 R W^2)\nonumber\\
&=&\alpha W^1RW^1+(1-\alpha)W^2RW^2\nonumber \\&& +\alpha (\alpha-1) \left[(W^1-W^2)R(W^1-W^2)\right]  \nonumber\\
&\preceq& \alpha W^1 R W^1+(1-\alpha)W^2 R W^2.
\end{IEEEeqnarray}
Therefore, $\mathscr{C}$ is convex.

For Statement (b), the following lemma is introduced first to relate a Riccati inequality condition to $\mathcal{H}_\infty$-norm attenuation.
\begin{lemma}~\cite{scherer1989h}\label{lemma:Riccati}
Given $\gamma >0$, if $[A_c, C_c]$ is observable, the closed-loop system is asymptotically stable and $\|H(s)\|_\infty \leq \gamma$ if and only if the Riccati inequality
\begin{IEEEeqnarray}{rl}
A_{c}^T P+P A_{c} + \gamma^{-2} P B_1 B_1^T P+ C_c^T C_c \preceq 0
\end{IEEEeqnarray}
has a symmetric positive definite solution $P=P^T\succ 0$.
\end{lemma}
\noindent \textbf{Proof of Lemma \ref{lemma:Riccati}:}
The proof is shown in~\cite{scherer1989h}.\hfill{\qed}

Notice that Assumption~\ref{assu:usual} implies that the pair $[A_c, C_c]$ is observable. Then, from Lemma \ref{lemma:Riccati}, there exists a symmetric positive definite solution $P=P^T\succ 0$ such that
\begin{IEEEeqnarray}{rl}
A_c^T P +P A_c    + \mu P B_1 B_1^T P+ C^T C+K^T D^T DK \preceq 0. \nonumber\\ \label{lemma:Riccatiexample}
\end{IEEEeqnarray}
Since $P$ is nonsingular, by pre-multiplying and post-multiplying $P^{-1}$ in \eqref{lemma:Riccatiexample}, we have
\begin{IEEEeqnarray}{rl}
P^{-1} A_c^T   +  A_c P^{-1}    + \mu   B_1 B_1^T + P^{-1} C^T C P^{-1} \nonumber\\+P^{-1} K^T D^T DK P^{-1} \preceq 0. \label{lemma:Riccatiexample2}
\end{IEEEeqnarray}
Denote $W_p=P^{-1}$, \eqref{lemma:Riccatiexample2} is equivalent to
\begin{IEEEeqnarray}{rCl}
A_c W_p+W_{p}A_c^T +W_{p}C^TCW_{p}&+W_{p}K^TD^TDKW_{p}   \nonumber\\&+\mu B_1B_1^T \preceq 0.
\end{IEEEeqnarray}
Meanwhile, from~\eqref{eq:theta_mu}, we have
\begin{IEEEeqnarray}{rCl}
\Theta_{1}(W,\mu)&=&AW_1-B_2W_2^T+W_1A^T-W_2B_2^T\nonumber\\&&+W_1C^TCW_1+W_2D^TDW_2^T +\mu B_1B_1^T. \IEEEeqnarraynumspace \label{eq:theta1}
\end{IEEEeqnarray}
Then, by setting $W_1=W_p$ and $W_2^T=KW_p$, we have $K=W_2^T W_1^{-1}$ and $\Theta_{1}(W,\mu)\preceq 0.$ 
It gives that $K=W_2^T W_1^{-1}$ is a feasible solution to ensure the stability with $\gamma$-attenuation~\cite{peres1994optimal}. By  substituting $K=W_2^T W_1^{-1}$ to~\eqref{eqn:partition}, we can construct \begin{IEEEeqnarray}{cl}W= \begin{bmatrix} W_1 & W_1 K^T\\KW_1 &W_3 \end{bmatrix}.\end{IEEEeqnarray}  By Schur complement, we can ensure $W \succeq 0$ by choosing $W_3\succeq K W_1 K^T$. Based on the analysis above, $K=W_2^T W_1^{-1}$ is a stabilizing gain generated from $(W, \mu) \in \mathscr{C}$, and it follows from Lemma \ref{lemma:Riccati} that $\|H(s)\|_\infty \leq \gamma$ is guaranteed~\cite{peres1994optimal}.

Statement (c) is direct consequence of Statement (b).  \hfill{\qed}

Then, it suffices to extend the above results to uncertain systems, and then we make the following assumption.

\begin{assumption}~\label{assu:uncer}
The parametric uncertainties are structural and convex-bounded.
\end{assumption}

Followed by Assumption~\ref{assu:uncer}, we have $F=\sum_{i=1}^N \xi_i F_i$,  $\xi_i \geq 0$, $\forall i=1,2, \cdots, N$,  and $\sum_{i=1}^N \xi_i=1$. Notice that $F$ belongs to a polyhedral domain, which can be expressed as a convex combination of the extreme matrices $F_i$, where \begin{IEEEeqnarray}{cl}
F_i=\begin{bmatrix}
A_i & -B_{2i}\\
0 & 0
\end{bmatrix} \in \mathbb{R}^{ p \times  p}.\end{IEEEeqnarray}  Then, define the matrical function in terms of each extreme vertice, where \begin{IEEEeqnarray}{cl}\Theta_i(W,\mu)=F_i W+WF_i^T+WRW+\mu Q, \end{IEEEeqnarray}  which can also be partitioned as
\begin{IEEEeqnarray}{cl} \Theta_i(W,\mu) = \begin{bmatrix} \Theta_{1i}(W,\mu) & \Theta_{2i}(W) \\
\Theta_{2i}^T(W) & \Theta_{3i}(W) \end{bmatrix},\end{IEEEeqnarray}
with $\Theta_{1i}(W,\mu) \in \mathbb{S}^{ n  },
\Theta_{2i}(W) \in \mathbb{R}^{ n \times m},
\Theta_{3i}(W) \in \mathbb{S}^{m  }$. Consequently,  a mapping between $W$ and $K$ can be constructed, and the results are shown in Theorem~\ref{thm:decentralizedconstraintrobust}.

\begin{theorem} \label{thm:decentralizedconstraintrobust}
Define the set $\mathscr{C}_U  = \{(W,\mu): W=W^T \succeq 0, \Theta_{1i}(W,\mu)\preceq 0, \mu > 0\}$. Then the following statements hold:
\begin{enumerate}[(a)]	
	\item Any $(W, \mu) \in \mathscr{C}_U$ generates a stabilizing gain $K=W_2^T W_1^{-1}$ that guarantees $\|H_i(s)\|_\infty \leq \gamma$, $\forall i=1,2, \cdots , N$, with $\gamma=1/\sqrt \mu>0$ under convex-bounded parametric uncertainties, where $\|H_i(s)\|_\infty$ represents the $\mathcal{H}_\infty$-norm with respect to
	the ith extreme system.
	\item At optimality, $(W^*,\mu^*)=\textup{argmax}\{\mu:(W, \mu)\in \mathscr{C}_U\}$ gives the optimal solution to the  $\mathcal{H}_\infty$ guaranteed cost control problem, with $K^*={W_2^*}^{T} {W_1^*}^{-1}$ and $\gamma^*=1/\sqrt{\mu^*}$.
\end{enumerate}
\end{theorem}

\noindent{\textbf{Proof of Theorem \ref{thm:decentralizedconstraintrobust}:}}
The proof is straightforward as it is an extension of Theorem \ref{thm:convex}, then it is omitted. \hfill{\qed}

\begin{remark}
Obviously $\gamma=1/\sqrt \mu$ is the upper bound to $\|{H}_i(s)\|_\infty$. For the uncertain systems, the upper bound is minimized at optimality; while for the precise systems, the upper bound is reduced to the optimal $\|{H}(s)\|_\infty$.
\end{remark}

\section{Symmetric Gauss-Seidel ADMM for $\mathcal{H}_\infty$ Guaranteed Cost Control}
\subsection{Formulation of the Optimization Problem}
Followed by the above analysis, the $\mathcal{H}_\infty$ guaranteed cost control problem can be formulated by the following convex optimization problem:
\begin{IEEEeqnarray}{rl}~\label{eq:opt1}
\displaystyle\operatorname*{maximize}_{(W,\mu)\in\mathbb S^{p}\times \mathbb R}
&\quad  \mu \nonumber\\
\operatorname*{subject\ to}&\quad
W \succeq 0\nonumber\\
&\quad \Theta_{1i}(W,\mu)   \preceq 0, \, \forall i=1,2,\cdots, N  \nonumber \\
&\quad\mu>0.
\end{IEEEeqnarray}

Define $V=\begin{bmatrix}
I_n & 0_{n \times m}
\end{bmatrix}$, and then \eqref{eq:opt1} can be equivalently expressed in the conventional form, where
\begin{IEEEeqnarray}{rl}~\label{eq:opt2}
\displaystyle\operatorname*{minimize}_{(W,\mu)\in\mathbb S^{p}\times \mathbb R}
&\quad -\mu \nonumber\\
\operatorname*{subject\ to}&\quad
W\in \mathbb S^{p}_+\nonumber\\
&\quad -V\left( F_i W+W  F_i^T+WRW+\mu Q\right)V^T  \nonumber\\    &\hspace{2.5cm} \in  \mathbb S^{n}_+, \, \forall i=1,2,\cdots,N\nonumber\\
&\quad\mu>0.
\end{IEEEeqnarray}
From Schur complement, for all $i=1,2,\cdots,N$, the second group of conic constraints in~\eqref{eq:opt2} can be equivalently expressed by
\begin{IEEEeqnarray}{rl}~\label{eq:schur}
\begin{bmatrix}
-VF_i WV^T-VW  F_i^TV^T-\mu VQV^T & VWR^{\frac{1}{2}}\\R^{\frac{1}{2}}WV^T & I_p
\end{bmatrix}  \succeq 0.  \nonumber\\
\end{IEEEeqnarray}
Then,~\eqref{eq:schur} can be further decomposed as
\begin{IEEEeqnarray}{rl}
\begin{bmatrix}
-VF_i\\R^{\frac{1}{2}}
\end{bmatrix}
W
\begin{bmatrix}
V^T & 0
\end{bmatrix}
+  \begin{bmatrix}V\\0\end{bmatrix}W \begin{bmatrix}-F_i^TV^T &R^{\frac{1}{2}} \end{bmatrix} \IEEEeqnarraynumspace\nonumber\\
+\mu  \begin{bmatrix}- VQV^T&0\\0&0\end{bmatrix}  +  \begin{bmatrix} 0 &0\\0&I_p \end{bmatrix}
\succeq 0.\IEEEeqnarraynumspace
\end{IEEEeqnarray}
For the sake of simplicity, we define
\begin{IEEEeqnarray}{rCl}
\mathcal G_i(W,\mu) =H_{i1}WH_{2}+H_{2}^TWH_{i1}^T+\mu H_{3}+H_{0},\IEEEeqnarraynumspace
\end{IEEEeqnarray}
where
\begin{gather}
H_{0}=\left[\begin{array}{cccc}0 &0\\0&I_p\end{array}\right]\in\mathbb{S}^{r}, \,
H_{i1}=\left[\begin{array}{cccc}-VF_i\\R^{\frac{1}{2}}\end{array}\right]\in\mathbb{R}^{r\times p}, \nonumber\\
H_{2}=\Big[V^T\quad 0\Big]\in\mathbb{R}^{p\times r}, \,
H_{3}=\left[\begin{array}{cccc}- VQV^T &0\\0&0\end{array}\right]\in\mathbb{S}^{r},\nonumber\\
\end{gather}
and  then the optimization problem is equivalently expressed as
\begin{IEEEeqnarray}{rl}
\displaystyle\operatorname*{minimize}_{(W,\mu)\in\mathbb S^{p}\times \mathbb R}
&\quad -\mu \nonumber\\
\operatorname*{subject\ to}&\quad
W\in \mathbb S^{p}_+\nonumber\\
&\quad	\mathcal G_i(W,\mu) \succeq 0, \, \forall i=1,2,\ldots,N \nonumber\\
&\quad \mu>0.
\end{IEEEeqnarray}

Then we introduce consensus variables $Y_0=W,\ Y_i=\mathcal G_i(W,\mu),\ \forall i=1,2,\cdots,N$, $Y_{N+1}=\mu$. Define a cone $\mathcal K$ as
\begin{IEEEeqnarray}{rCl}
\mathcal K=\mathbb S^{p}_+\times \underbrace{\mathbb S^{r}_+\times \mathbb S^{r}_+\times \cdots\times \mathbb S^{r}_+}_N\times \mathbb R_+,
\end{IEEEeqnarray}
and also the corresponding linear space $\mathcal X$ as
\begin{IEEEeqnarray}{rCl}
\mathcal X=\mathbb S^{p}\times \underbrace{\mathbb S^{r}\times \mathbb S^{r}\times \cdots\times \mathbb S^{r}}_N\times \mathbb R.
\end{IEEEeqnarray}
Notably, since the positive semi-definite cone is self-dual, it follows that $\mathcal K=\mathcal K^*\subset \mathcal X$, where $\mathcal{K}^*$ represents the dual of $\mathcal{K}$. Besides, define a linear mapping $\mathcal H:\mathbb S^p\times\mathbb R\rightarrow \mathcal X$, where
$
\mathcal H(W,\mu)=\big(W,\mathcal G_1(W,\mu),\mathcal G_2(W,\mu),\dotsm,\mathcal G_N(W,\mu),\mu\big),
$
and define the corresponding vector $Y=(Y_0,Y_1,\dotsm,Y_{N+1})$ in the given space $\mathcal X$. Then the optimization problem can be transformed into the following compact form:
\begin{IEEEeqnarray}{rl}
\displaystyle\operatorname*{minimize}_{(W,\mu)\in\mathbb S^{p}\times \mathbb R}
&\quad -\mu+\delta_{\mathcal K}(Y)\nonumber\\
\operatorname*{subject\ to}
&\quad Y-\mathcal H(W,\mu)=0,
\end{IEEEeqnarray}
where $\delta_\mathcal K(Y)$ is the indicator function in terms of the convex cone $\mathcal K$, which is given by
\begin{IEEEeqnarray}{rl}
\delta_\mathcal K(Y)=\begin{cases}
0 & \text{if } Y \in\mathcal K  \\
+\infty  & \text{otherwise.}
\end{cases}
\end{IEEEeqnarray}

In order to deal with the problems leading to the large-scale optimization, a serial computation technique is introduced. Before we present the optimization procedures in detail, define the augmented Lagrangian as
\begin{IEEEeqnarray}{rCl}
\mathcal L_\sigma(Y,W,\mu;Z) &=& -\mu+\delta_{\mathcal K}(Y)
+\frac{\sigma}{2}\big\|Y-\mathcal H(W,\mu)\nonumber\\&&
+\sigma^{-1}Z\big\|^2-\frac{1}{2\sigma}\|Z\|^2,
\end{IEEEeqnarray}
where $Z=(Z_0,Z_1,\dotsm,Z_{N+1})\in\mathcal K^*$ is the vector of the Lagrange multipliers.

\subsection{Symmetric Gauss-Seidel ADMM Algorithm}
The numerical procedures of the symmetric Gauss-Seidel ADMM algorithm are given below, where $Y$, $W$, $\mu$, and $Z$ are updated through an iterative framework. By using the proposed algorithm, $Y$ can be updated by parallel computation, such that high efficiency and feasibility can be ensured even for the large-scale optimization problems, $W$ and $\mu$ are updated in a serial framework such that there is an explicit solution to the sub-problem in terms of each one of both variables, and finally the Lagrange multiplier $Z$ is updated.

\subsection*{Step 1. Initialization}
For initialization, the following parameters and matrices need to be selected first:  $\tau=1.618$, in fact, $\tau$ can be chosen within $(0,(1+\sqrt 5)/2)$; $\sigma$ is chosen as a positive real number;  $\left(Y^0;W^0,\mu^0\right)\in\mathcal X\times \mathbb S^{p} \times \mathbb R$ and $Z^0\in\mathcal X$; $\epsilon>0$. Then, set the iteration index $k=0$.

\subsection*{Step 2. Update of $Y$}
In the following text, we define $\partial (\cdot)$ as the sub-differential operator. Also, we define $\partial_\mathcal{S} (\cdot)$ as the sub-differential operator with respect to the variable ${S}$.
Since the sub-problem for updating the variable $Y$ is unconstrained, the optimality condition is given by
\begin{IEEEeqnarray}{rCl}
0\in \partial_Y\mathcal L_\sigma\Big(Y,W^k,\mu^k;Z^k\Big).
\end{IEEEeqnarray}
Notice that when there are a large number of uncertainties in the given system, the sub-problem is not easy to solve in terms of the whole variable vector $Y$. Therefore, a parallel computation technique is proposed to cater to this practical constraint. To solve this problem with the parallel computation technique, we rearrange the augmented Lagrangian into the following form:
\begin{IEEEeqnarray}{rCl}
&&\quad \mathcal L_\sigma(Y,W,\mu;Z)\nonumber\\ &=& -\mu+\delta_{\mathbb S_+^p}(Y_{0})+
\sum_{i=1}^N \delta_{\mathbb S_+^{r}}(Y_{i})+\delta_{\mathbb R_+}(Y_{N+1})\nonumber\\
&&+\frac{\sigma}{2}\|Y_0-W+\sigma^{-1}Z_0\|^2
+\sum_{i=1}^N  \frac{\sigma}{2}\|Y_i-\mathcal G_i(W,\mu)\nonumber\\
&&+\sigma^{-1}Z_i\|^2+\frac{\sigma}{2}\|Y_{N+1}-\mu+\sigma^{-1}Z_{N+1}\|^2
-\frac{1}{2\sigma}\|Z\|^2, \nonumber\\
\end{IEEEeqnarray}
where $\delta_{\mathbb S_+^p}(\cdot),\ \delta_{\mathbb S_+^n}(\cdot)$, and $\delta_{\mathbb R_+}(\cdot)$ are the indicator functions in terms of the $p$ dimensional positive semi-definite cone, $n$ dimensional positive semi-definite cone, and positive cone of the real numbers, respectively.

First of all, we consider the optimality condition to the sub-problem in terms of the variable $Y_{N+1}$, which is given by
\begin{IEEEeqnarray}{rCl}
0&\in& \partial_{Y_{N+1}}\mathcal L_\sigma\Big(Y,W^k,\mu^k;Z^k\Big)\nonumber\\
&\in&\partial\delta_{\mathbb R_+}(Y_{N+1})+\sigma\Big(Y_{N+1}-\mu^k+\sigma^{-1}Z^k_{N+1}\Big).\IEEEeqnarraynumspace
\end{IEEEeqnarray}

To determine the projection operator $\Pi_\mathcal C(\cdot)$ with respect to the convex cone $\mathcal C$, the following theorem is used.
\begin{theorem}~\cite{ma2021optimal}~\label{thm:proj}
The projection operator $\Pi_\mathcal C(\cdot)$ with respect to the convex cone $\mathcal C$ can be expressed as
\begin{IEEEeqnarray}{rCl}
\Pi_{\mathcal C}=(I+\alpha \partial \delta_{\mathcal C})^{-1},
\end{IEEEeqnarray}
where $\alpha \in \mathbb R$ can be an arbitrary real number.
\end{theorem}
\noindent \textbf{Proof of Theorem \ref{thm:proj}:}
The complete proof can be found in~\cite{ma2021optimal}. \hfill{\qed}

Therefore, we have
\begin{IEEEeqnarray}{rCl}\label{equation: ADMM1_1}
\mu^k-\sigma^{-1}Z^k_{N+1}&\in& (\sigma^{-1}\partial\delta_{\mathbb R_+}+I)(Y_{N+1})\nonumber\\
Y_{N+1}^{k+1}&=&\Pi_{\mathbb R_+}\Big(\mu^k-\sigma^{-1}Z^k_{N+1}\Big).
\end{IEEEeqnarray}

To calculate the projection operator in terms of the positive semi-definite convex cone explicitly, the following lemma is introduced.

\begin{lemma}~\cite{ma2021optimal}\label{lemma:proj}
Projection onto the positive semi-definite cone can be computed explicitly. Let $X = \sum_{i=1}^n \lambda_i v_i v_i^T\in\mathbb S^n$ be the eigenvalue decomposition of the matrix $X$ with the eigenvalues satisfying $\lambda_1 \geq \lambda_2 \geq \dots \geq \lambda_n$, where $v_i$ denotes the eigenvector corresponding to the $i$th eigenvalue. Then the projection onto the positive semi-definite cone of the matrix $X$ can be expressed by
\begin{IEEEeqnarray}{rCl}
\Pi_{\mathbb S^n_+} (X) = \sum_{i=1}^n \max\left\{\lambda_i, 0\right\} v_i v_i^T.
\end{IEEEeqnarray}
\end{lemma}
\noindent \textbf{Proof of Lemma \ref{lemma:proj}:}
The proof is shown in~\cite{ma2021optimal}.
\hfill{\qed}

Then we consider the optimality condition to the sub-problem in terms of the variable $Y_i$, $\forall i=N,N-1,\cdots,1$, where
\begin{IEEEeqnarray}{rCl}
0&\in& \partial_{Y_{i}}\mathcal L_\sigma\Big(Y,W^k,\mu^k;Z^k\Big)\nonumber\\
&\in&\partial\delta_{\mathbb R_+^{r}}(Y_i)+\sigma\Big(Y_i-\mathcal G_i\left(W^k,\mu^k\right)+\sigma^{-1}Z^k_{i}\Big).\IEEEeqnarraynumspace
\end{IEEEeqnarray}
Therefore, we have
\begin{IEEEeqnarray}{rCl}\label{equation: ADMM1_2}
\mathcal G_i(W^k,\mu^k)-\sigma^{-1}Z^k_{i}&\in& \left(\sigma^{-1}\partial\delta_{\mathbb R_+^{r}}+I\right)(Y_i)\nonumber\\
Y_{i}^{k+1}&=&\Pi_{\mathbb S_+^{r}}\Big(\mathcal G_i\left(W^k,\mu^k\right)-\sigma^{-1}Z^k_{i}\Big).\nonumber\\
\end{IEEEeqnarray}
Then we consider the optimality condition to the sub-problem in terms of the variable $Y_0$, where
\begin{IEEEeqnarray}{rCl}
0&\in& \partial_{Y_{i}}\mathcal L_\sigma\Big(Y,W^k,\mu^k;Z^k\Big)\nonumber\\
&\in&\partial\delta_{\mathbb R_+^{p}}(Y_0)+\sigma\Big(Y_0-W^k+\sigma^{-1}Z^k_{0}\Big).
\end{IEEEeqnarray}
Therefore, we have
\begin{IEEEeqnarray}{rCl}\label{equation: ADMM1_3}
W^k-\sigma^{-1}Z^k_{0}&\in& (\sigma^{-1}\partial\delta_{\mathbb R_+^{p}}+I)(Y_0)\nonumber\\
Y_{0}^{k+1}&=&\Pi_{\mathbb S_+^{p}}\Big(W^k-\sigma^{-1}Z^k_{0}\Big).
\end{IEEEeqnarray}
\begin{remark}
Notice that each projection can be computed independently, which means that no more information is required to obtain the projection of each variable onto the corresponding convex cone, except for the value of the same variable in the last iteration. Therefore, the projection of $Y$ onto the convex cone $\mathcal K$ can be obtained by solving a group of separate sub-problems.
\end{remark}

\subsection*{Step 3. Update of $W$ and $\mu$}
The optimality conditions in terms of the sub-problem of the variable set $(W,\mu)$ are given by
\begin{IEEEeqnarray}{rCl}
\left\{\begin{array}{l}
0\in\partial_{W} \mathcal L_\sigma(Y,W,\mu;Z)\\
0\in\partial_{\mu} \mathcal L_\sigma(Y,W,\mu;Z).
\end{array}\right.
\end{IEEEeqnarray}
To solve this sub-problem efficiently, the symmetric Gauss-Seidel technique is introduced. Before the optimality condition is given,  the following lemma is presented which determines the derivation of a norm function with a specific structure.
\begin{lemma}~\label{lemma:norm_deri}
Given a norm function
\begin{IEEEeqnarray}{rCl}\label{eqn:normmm}
\mathcal H_i(W)=\| H_i(W)\|^2,
\end{IEEEeqnarray}
where
\begin{IEEEeqnarray}{rCl}
H_i(W)=H_{i1}WH_{2}+H_{2}^TWH_{i1}^T+\mu H_{3}+H_{0},
\end{IEEEeqnarray}	
$H_{0},H_{i1}$, $H_2$, and $H_3$ are given matrices with appropriate dimensions. Then it follows that
\begin{IEEEeqnarray}{rCl}
\frac{\partial \mathcal H_i(W)}{\partial W}=2H_{i1}^T H_i(W)H_2^T+2H_2H_i(W)H_{i1}.
\end{IEEEeqnarray}
\end{lemma}
\noindent \textbf{Proof of Lemma \ref{lemma:norm_deri}:}
The derivative of the matrix norm function in the form of~\eqref{eqn:normmm} can be obtained by using some properties of derivative of trace operator. The procedures are simple but tedious, so the proof is omitted.
\hfill{\qed}

On the basis of the symmetric Gauss-Seidel technique, the optimality conditions to the sub-problems in the backward sweep and the forward sweep are given in Step 3.1 and Step 3.2, respectively.

\subsection*{Step 3.1. Symmetric Gauss-Seidel Backward Sweep}
From Lemma~\ref{lemma:norm_deri}, we can easily obtain the derivatives of the norm functions with respect to the corresponding variables.
Consider the optimality condition of the sub-problem in terms of the variable $\mu$, it follows that
\begin{IEEEeqnarray}{rCl}
0&\in&\partial_{\mu} \mathcal L_\sigma(Y^{k+1},W^k,\mu;Z^k)\nonumber\\
&=&-1+\sigma\sum_{i=1}^{N}\Big\langle \mathcal G_i(W^k,\mu)-Y_i^{k+1}-\sigma^{-1}Z_i^{k},H_3\Big\rangle \nonumber\\
&&+\sigma\left(\mu-Y_{N+1}^{k+1}-\sigma^{-1}Z_{N+1}^k\right)\nonumber\\
&=&-1+\sigma\sum_{i=1}^{N}\Big\langle
H_{i1}W^k H_2+H_2^TW^kH_{i1}^T-Y_i^{k+1}\nonumber\\
&&-\sigma^{-1}Z_i^{k},H_3\Big\rangle+\sigma\left(-Y_{N+1}^{k+1}-\sigma^{-1}Z_{N+1}^k\right)\nonumber\\
&&+\sigma N\langle H_0,H_3\rangle +\mu \sigma \Big[N \operatorname{Tr} \big(H_3^2\big)+1\Big].
\end{IEEEeqnarray}
Therefore, we have
\begin{IEEEeqnarray}{rCl}\label{equation: ADMM2_1}
\bar \mu^{k+1}&=&\sigma^{-1}\Big[N \operatorname{Tr} \big(H_3^2\big)+1\Big]^{-1}\Bigg(
1-\sigma\sum_{i=1}^{N}\Big\langle H_{i1}W^k H_2 \nonumber\\
&&+H_2^TW^kH_{i1}^T-Y_i^{k+1}-\sigma^{-1}Z_i^{k},H_3\Big\rangle
\nonumber\\&&+\sigma Y_{N+1}^{k+1}+Z_{N+1}^k
-\sigma N\langle H_0,H_3\rangle \Bigg).
\end{IEEEeqnarray}
Then we consider the optimality condition of the sub-problem in terms of the variable $W$, we have
\begin{IEEEeqnarray}{rCl}
0&\in&\partial _W\mathcal L_\sigma (Y^{k+1},W,\bar \mu^{k+1};Z^k)\nonumber\\
&=&\sigma(W-Y_0-\sigma^{-1}Z_0)\nonumber\\
&&+\sigma\sum_{i=1}^N\Bigg[ H_{i1}^T\Big(\mathcal G_i(W,\bar \mu^{k+1})-Y_i^{k+1}-\sigma^{-1}Z_i^{k}\Big)H_2^T\nonumber\\
&&+H_2\Big(\mathcal G_i(W,\bar \mu^{k+1})-Y_i^{k+1}-\sigma^{-1}Z_i^{k}\Big)H_{i1}\Bigg]\nonumber\\
&=&W-Y_0-\sigma^{-1}Z_0\nonumber\\
&&+\sum_{i=1}^N\Bigg[ H_{i1}^T\Big(\bar \mu^{k+1} H_{3}+H_{0}-Y_i^{k+1}-\sigma^{-1}Z_i^{k}\Big)H_2^T\nonumber\\
&&+H_2\Big(\bar \mu^{k+1} H_{3}+H_{0}-Y_i^{k+1}-\sigma^{-1}Z_i^{k}\Big)H_{i1}\Bigg]\nonumber\\
&&+\sum_{i=1}^N\Bigg[H_{i1}^T\Big(H_{i1}WH_{2}+H_{2}^TWH_{i1}^T\Big)H_2^T \nonumber\\
&&+H_2\Big(H_{i1}WH_{2}+H_{2}^TWH_{i1}^T\Big)H_{i1}\Bigg].
\end{IEEEeqnarray}
To obtain $W$ explicitly, the vectorization technique is utilized, then define
\begin{IEEEeqnarray}{rCl}
T_0 &=&-Y_0-\sigma^{-1}Z_0\nonumber\\
&&+\sum_{i=1}^N\Bigg[ H_{i1}^T\Big(\bar \mu^{k+1} H_{3}+H_{0}-Y_i^{k+1}-\sigma^{-1}Z_i^{k}\Big)H_2^T\nonumber\\
&&+H_2\Big(\bar \mu^{k+1} H_{3}+H_{0}-Y_i^{k+1}-\sigma^{-1}Z_i^{k}\Big)H_{i1}\Bigg],
\end{IEEEeqnarray}
and then it follows that
\begin{IEEEeqnarray}{rCl}~\label{eq:vec}
0 &=& T_0 +W+\sum_{i=1}^N\Bigg[H_{i1}^T\Big(H_{i1}WH_{2}+H_{2}^TWH_{i1}^T\Big)H_2^T \nonumber\\
&&+H_2\Big(H_{i1}WH_{2}+H_{2}^TWH_{i1}^T\Big)H_{i1}\Bigg].
\end{IEEEeqnarray}
It is straightforward that \eqref{eq:vec} is equivalent to
\begin{IEEEeqnarray}{rCl}
0&=&\operatorname{vec}(T_0)+\Bigg[I+\sum_{i=1}^{N}\Big[
(H_2H_2^T)\otimes(H_{i1}^TH_{i1}) \nonumber\\&&+(H_2 H_{i1})\otimes (H_{i1}^TH_2^T)+(H_{i1}^TH_2^T)\otimes(H_2H_{i1})\nonumber\\&&+(H_{i1}^TH_{i1})\otimes (H_2H_2^T)
\Big]\Bigg]\operatorname{vec}(W).
\end{IEEEeqnarray}
Then it follows that
\begin{IEEEeqnarray}{l}\label{equation: ADMM2_2}
\operatorname{vec}(W^{k+1})\nonumber\\
=-\Bigg[I+\sum_{i=1}^{N}\Big[
(H_2H_2^T)\otimes(H_{i1}^TH_{i1}) \nonumber\\\quad+(H_2 H_{i1})\otimes (H_{i1}^TH_2^T)+(H_{i1}^TH_2^T)\otimes(H_2H_{i1})\nonumber\\\quad+(H_{i1}^TH_{i1})\otimes (H_2H_2^T)
\Big]\Bigg]^{-1}\operatorname{vec}(T_0).
\end{IEEEeqnarray}
In this way, $W^{k+1}$ can be obtained by performing the inverse vectorization.

\subsection*{Step 3.2. Symmetric Gauss-Seidel Forward Sweep}
\begin{IEEEeqnarray}{rCl}\label{equantion: ADMM2_3}
\mu^{k+1}&=&\sigma^{-1}\Big[N \operatorname{Tr} \big(H_3^2\big)+1\Big]^{-1}\Bigg(
1-\sigma\sum_{i=1}^{N}\Big\langle
H_{i1}W^{k+1} H_2\nonumber\\&&+H_2^TW^{k+1}H_{i1}^T-Y_i^{k+1}-\sigma^{-1}Z_i^{k},H_3\Big\rangle
\nonumber\\&&+\sigma Y_{N+1}^{k+1}+Z_{N+1}^k
-\sigma N\langle H_0,H_3\rangle
\Bigg).
\end{IEEEeqnarray}
\begin{remark}
By using the symmetric Gauss-Seidel technique, the optimization procedures for the variable $W$ and the variable $v$ can be separated. The computational complexity is reduced significantly, because no matrical equation is required to be solved with the proposed algorithm comparing with the conventional ADMM counterpart.
\end{remark}

\subsection*{Step 4. Update of $Z$.}
\begin{IEEEeqnarray}{rCl}\label{equation: ADMM3}
Z^{k+1}=Z^k+\tau\sigma\Big(Y^{k+1}-\mathcal H\left(W^{k+1},\mu^{k+1}\right)\Big).
\end{IEEEeqnarray}

\subsection*{Step 5. Check the Stopping Criterion}

To derive the stopping criterion for the numerical procedures,  define the Lagrangian as
\begin{IEEEeqnarray}{rCl}
\mathcal L(W,\mu,Y;Z) = -\mu+\delta_{\mathcal K}(Y)+\langle Z,Y-\mathcal H(W,\mu)\rangle,  \nonumber\\
\end{IEEEeqnarray}
and then the KKT optimality conditions are given by
\begin{IEEEeqnarray}{rCl}\label{eq:KKT1}
\left\{\begin{array}{l}
0\in\partial _W\mathcal L(W,\mu,Y;Z)\\
0\in\partial _\mu\mathcal L(W,\mu,Y;Z)\\
0\in\partial _Y\mathcal L(W,\mu,Y;Z)\\
Y-\mathcal H(W,\mu)=0.
\end{array}\right.
\end{IEEEeqnarray}
It is straightforward that the relative residual errors are given by
\begin{IEEEeqnarray*}{rCl}
\operatorname{err}^k_W&=&\dfrac{\left\|Z_0^k+\sum_{i=1}^N\Big( H_{i1}^TZ_i^kH_2^T+H_{2}Z_i^kH_{i1}\Big) \right\|}{1+\left\|Z_0^k\right\|+\sum_{i=1}^M\left\| H_{i1}^TZ_iH_2^T+H_{2}Z_iH_{i1}\right\|}\\
\operatorname{err}^k_\mu&=&\dfrac{\left\|1+Z_{N+1}^k+\operatorname{Tr}\Big( \sum_{i=1}^NZ_i^kH_3\Big) \right\|}{2}\\
\operatorname{err}^k_Y&=&\dfrac{\left\|Y^k-\Pi_{\mathcal K}(Y^k-Z^k)\right\|}{1+\| Y^k\| + \|Z^k\|}\\
\operatorname{err}^k_{eq}&=&\dfrac{\left\|Y^k-\mathcal H(W^k,\mu^k) \right\|}{1+\|Y^k\|+\left\| \mathcal H(W^k,\mu^k)\right\|}.\yesnumber
\end{IEEEeqnarray*}
Define the relative residual error as
\begin{IEEEeqnarray*}{rCl}\label{equation: rre}
\operatorname{err}^k&=&\max\left\{\operatorname{err}_W^k,\operatorname{err}_\mu^k,\operatorname{err}_Y^k,\operatorname{err}_{eq}^k\right\}.\yesnumber
\end{IEEEeqnarray*}

According to the KKT optimality conditions, when the optimization variables are approaching their optimums, the relative residual errors are approaching zero. However, because of the numerical errors, the relative residual errors converge to a very small number instead of zero. Therefore, a small number $\epsilon$ is chosen as the stopping criterion, and when the stopping criterion $\operatorname{err}^k<\epsilon$ is satisfied, the current variables are at optimality.

\begin{remark}
The precision of the optimality can be increased with a tightened stopping criterion, though it would sacrifice the computational efficiency. 
\end{remark}

To this point, these numerical procedures are summarized by Algorithm~1.
\begin{algorithm}\label{algo}
\caption{Symmetric Gauss-Seidel ADMM for $\mathcal{H}_\infty$ guaranteed cost control}
\begin{algorithmic}[1]\label{algorithm: back_tracking_line_search}
	\REQUIRE
	Initialize the parameters $\sigma$, $\tau$, and $\epsilon$, the matrices $\left(Y^0;W^0,\mu^0\right)$ and $Z^0$. Set the iteration index $k=0$. For $k=0,1,2,\cdots$, perform the $k$th iteration
	\ENSURE $K^*$, $\gamma^*$
	\WHILE {\TRUE}
	\STATE   Determine $Y^{k+1}$ by~\eqref{equation: ADMM1_1}, \eqref{equation: ADMM1_2}, and \eqref{equation: ADMM1_3}.
	\STATE   Determine $\bar \mu^{k+1}$ and $\operatorname{vec}(W^{k+1})$ by~\eqref{equation: ADMM2_1} and~\eqref{equation: ADMM2_2}, respectively, and do the inverse vectorization to $\operatorname{vec}(W^{k+1})$ such that $W^{k+1}$ can be determined.
	\STATE   Determine $\mu^{k+1}$ by \eqref{equantion: ADMM2_3}.
	\STATE   Determine $z^{k+1}$ by \eqref{equation: ADMM3}.
	\STATE   Determine $\operatorname{err}^{k+1}$ by \eqref{equation: rre}.
	\IF {$\operatorname{err}^{k+1}<\epsilon$}
	\STATE {$K^*=(W_2^{k+1})^T(W_1^{k+1})^{-1}$}
	\STATE {$\gamma^*=1/\sqrt{\mu^{k+1}}$}
	\STATE \textbf{break}
	\ENDIF
	\ENDWHILE
	\RETURN  $K^*$, $\gamma^*$
\end{algorithmic}
\end{algorithm}

\subsection{Convergence Analysis and Computational Burden }
It is well-known that the conventional ADMM algorithm with a two-block structure can converge to the optimum linearly under mild assumptions~\cite{bertsekas1989parallel}. However, for the directly extended ADMM optimization with a multi-block structure, even with very small step size, the convergence cannot be ensured for particular optimization problems~\cite{chen2016direct}. To overcome this limitation, the symmetric Gauss-Seidel algorithm is proposed, and it can be proved that a linear convergence rate is guaranteed under the assumptions in terms of the linear-quadratic non-smooth cost function, such that the practicability and efficiency of ADMM technique to solve the large-scale optimization problems is significantly improved. Since the linear non-smooth cost function is a special case of the linear-quadratic non-smooth cost function, it is straightforward that the convergence of the proposed algorithm is guaranteed. Additionally, matrix operations are adopted in each iteration, which facilitate the use of vectorized implementation and reduce the computational burden significantly. 

\subsection{Discussion}
The methodology presented in this work can be broadly used when the controller gain is under prescribed structural constraints. For example, the synthesis of a decentralized controller can be determined by relating the decentralized structure to the certain equality constraints in the parameter space. Also, any controller with sparsity constraints can be converted to the decentralized constraints by factorization~\cite{ma2019parameter}. Further extensions also include the output feedback problem, which can be reformulated as a state feedback problem with a structural constraint~\cite{geromel1993convex}. These constraints can be simply added to the optimization problem to be solved by the symmetric Gauss-Seidel ADMM algorithm.


\section{Illustrative Examples}
To illustrate the effectiveness of the above results, two examples are presented. Example~\ref{exam:1} is reproduced from~\cite{petersen1987procedure}, which presents a robust controller design problem for an F4E fighter aircraft with a precise model in the longitudinal short period mode. Example~\ref{exam:2} presents a controller design problem with parametric uncertainties, which leads to a large-scale optimization problem. In this example, the state matrix and the control input matrix are randomly chosen such that their elements are stochastic variables uniformly distributed over $[0, 1]$, and parametric uncertainties with a variation of $\pm$20\% are applied to all parameters in the state matrix and the control input matrix.

In these examples, the optimization algorithm is implemented in Python 3.7.5 with Numpy 1.16.4, and executed on a computer with 16GB RAM and a 2.2GHz i7-8750H processor (6 cores).
For Example~\ref{exam:1}, the parameters for initialization is given by: $\sigma=0.001$, $\tau=0.618$, $\epsilon = 10^{-4}$, $Y^0=0$, $W^0=0$, $\mu^0=0$, $Z^0=0$; for Example~\ref{exam:2}, $\sigma$ and $\epsilon$ are set as 0.1 and $5\times10^{-4}$ for better convergence with all the other parameters remaining the same.

\begin{example}\label{exam:1}
Denote $x=	\begin{bmatrix}
N_z&  q&  \delta_e
\end{bmatrix}^T$, where $N_z$, $q$, and $\delta_e$ represent the normal acceleration, pitch rate, and elevation angle, respectively, and then the state space model of the aircraft is given by
\begin{IEEEeqnarray*}{rl}
\dot x &=A {x} +B_2 u+B_1 w\\
z&=Cx+Du\\
u&=-Kx,\nonumber
\end{IEEEeqnarray*}
where
\begin{gather*}\nonumber
A =
\begin{bmatrix}
-0.9896&  17.41& 96.15   \\
0.2648 & -0.8512 &-11.39   \\
0 &0& -30
\end{bmatrix},  \,
B_2  =
\begin{bmatrix}
-97.78\\ 0\\ 30
\end{bmatrix}, 
\\
B_1= \begin{bmatrix}
1 & 0 & 0\\ 0 & 1 & 0\\ 0 & 0 & 1
\end{bmatrix}, \,
C=\begin{bmatrix}
1 & 0 & 0\\ 0 & 1 & 0\\ 0 & 0 & 0
\end{bmatrix}, \,
D=\begin{bmatrix}
0 \\0\\1
\end{bmatrix}.\nonumber
\end{gather*}
\end{example}

\begin{figure}[t!]
\centering
\includegraphics[trim=0 0 0 0,width=1\columnwidth]{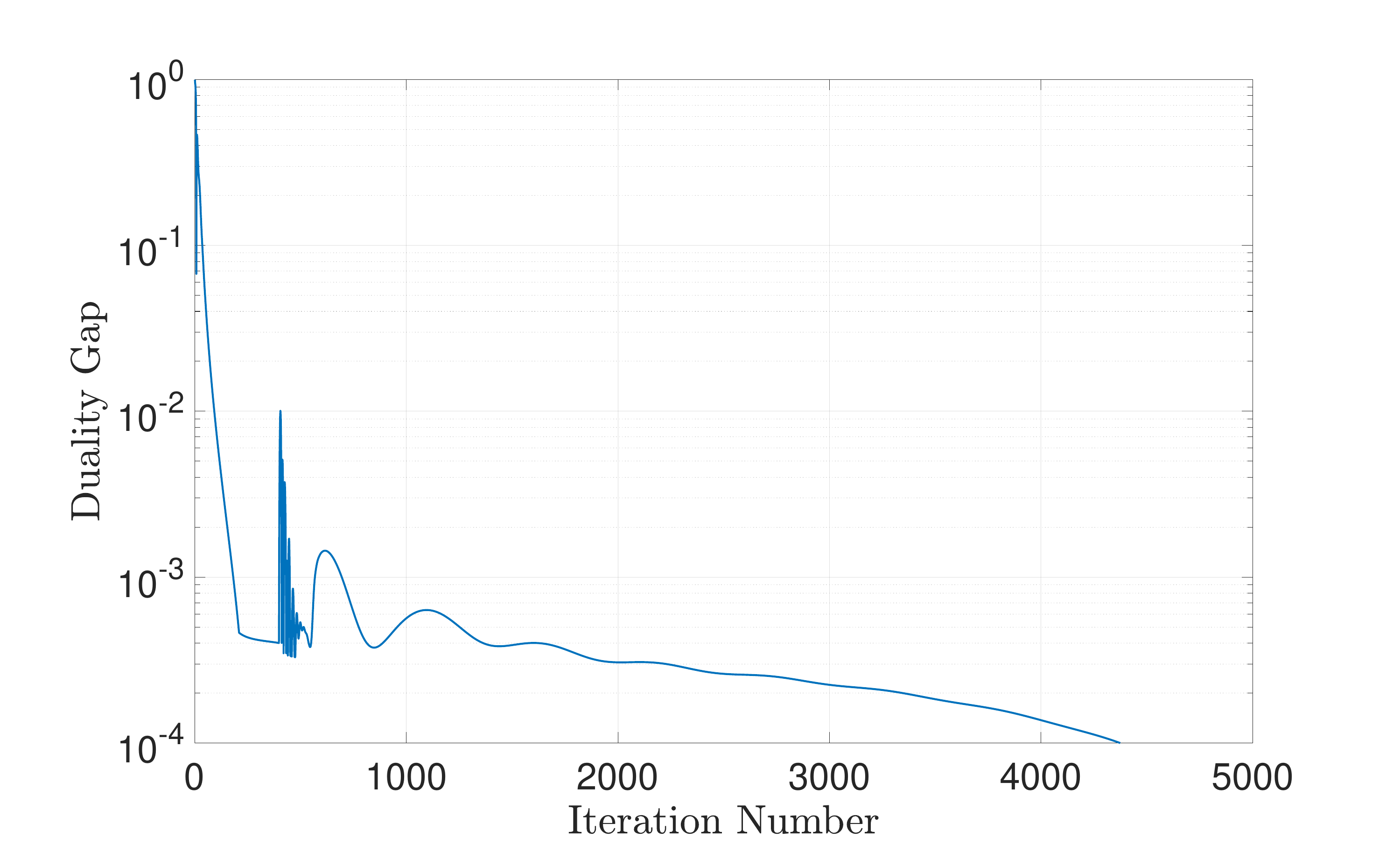}
\caption{Duality gap during iterations in Example 1}
\label{fig:cost1}
\end{figure}

In this example, we performed 5 trials using the proposed approach and recorded the total computation time, which is given by 4.3935 s, 4.6554 s, 4.1892 s, 4.7928 s, and 4.2933 s. In this case, the average computation time is 4.4648 s. Also, the change of the duality gap is shown in Fig.~\ref{fig:cost1}. At optimality, $W^*$ and $\mu^*$ are obtained, where
\begin{IEEEeqnarray*}{rCl}
W^* &=& \begin{bmatrix}
41.2179 & -3.9386 & -12.5019 & 4.7155\\
-3.9386 & 0.8802 & 0.7774 & -0.8563\\
-12.5019 & 0.7774 & 4.2692 & -1.6152\\
4.7155 & -0.8563 & -1.6152 & 127.8537
\end{bmatrix},
\end{IEEEeqnarray*}	and
\begin{IEEEeqnarray*}{rCl}
\mu^* &=& 4.4342.
\end{IEEEeqnarray*}	
It can be verified that all the constraints in the optimization problem are exactly satisfied. Then, the optimal controller gain is given by
\begin{IEEEeqnarray*}{l}
K^*=	 \begin{bmatrix}
-1.4754 & -4.0811 & -3.9557
\end{bmatrix},
\end{IEEEeqnarray*}
and the minimum level of disturbance attenuation is given by
\begin{IEEEeqnarray*}{l}
\gamma^*=0.4749.
\end{IEEEeqnarray*}

\begin{figure}[t!]
\centering
\includegraphics[trim=0 0 0 0,width=1\columnwidth]{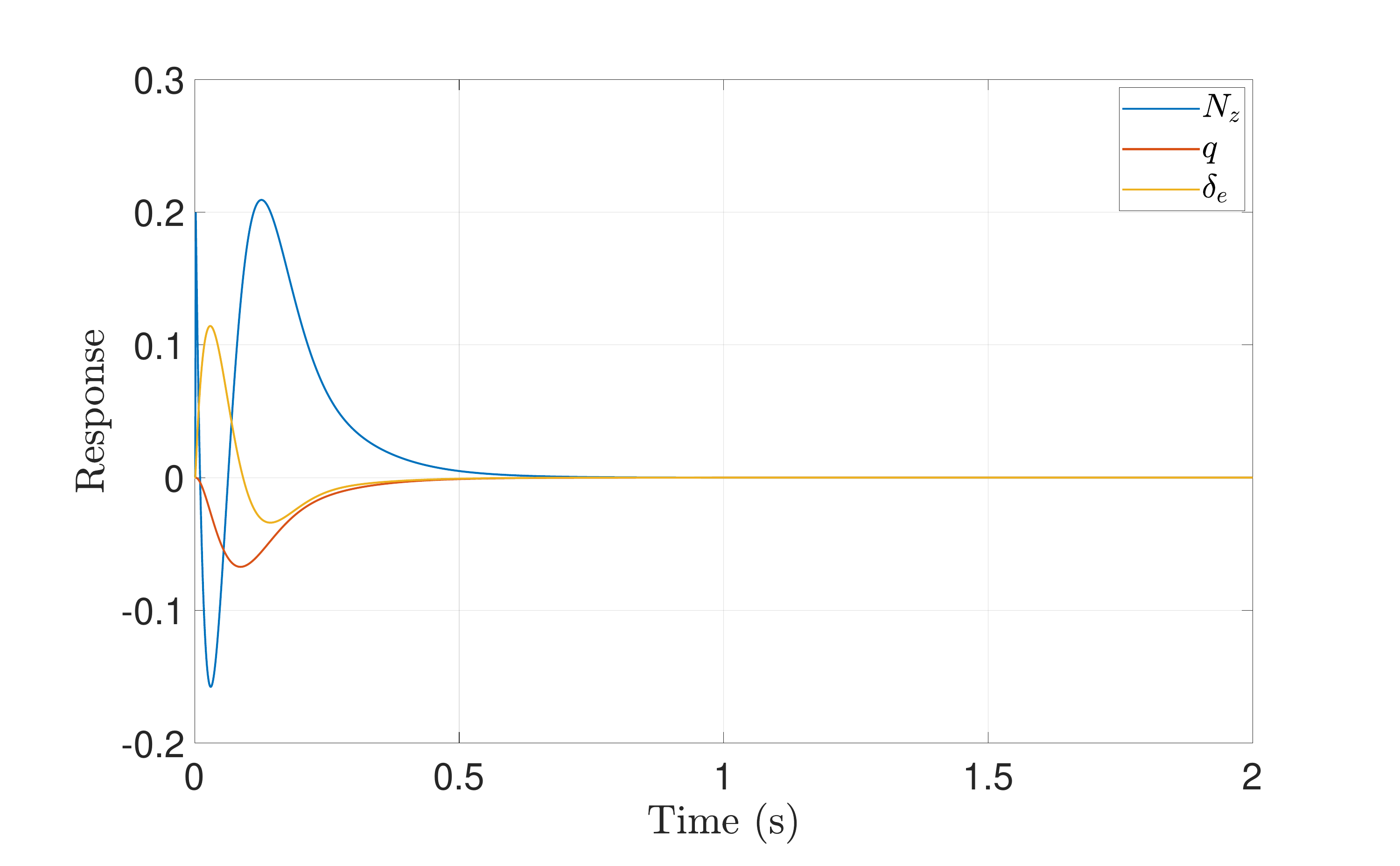}
\caption{System response in Example 1}
\label{fig:response1}
\end{figure}

\begin{figure}[t!]
\centering
\includegraphics[trim=0 0 0 0,width=1\columnwidth]{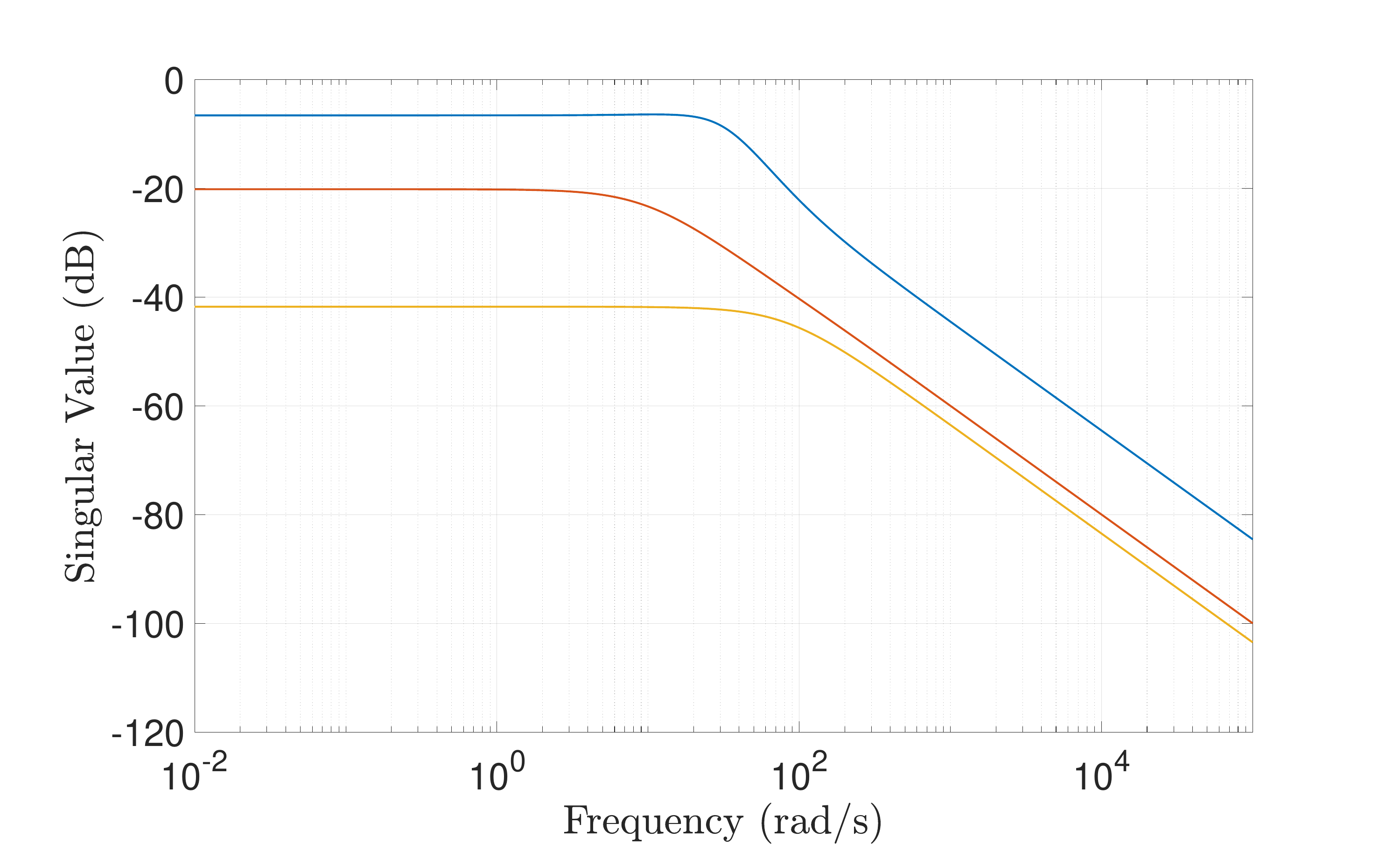}
\caption{Singular value diagram in Example 1}
\label{fig:sv1}
\end{figure}

In the simulation, consider $w$ as a vector of the impulse disturbance, and then the responses of normal acceleration, pitch rate, and elevation angle are shown in Fig.~\ref{fig:response1}. It can be easily verified that the closed-loop stability is suitably ensured. Beside, the singular value diagram of $H(s)$ is shown in Fig.~\ref{fig:sv1}. In the diagram, the maximum singular value is given by -6.43 dB, which is equivalent to 0.4770 in magnitude. It is almost the same as $\gamma^*$ that we have computed, and this is tallied with the condition that there is no parametric uncertainty in the model. 


\begin{example}\label{exam:2}
Consider $x=\begin{bmatrix}
x_1 & x_2 
\end{bmatrix}^T$ and a linear system
\begin{IEEEeqnarray*}{rl}
\dot {{x}}&=A {x} +B_2 u+B_1 w\\
z&=Cx+Du\\
u&=-Kx,\nonumber
\end{IEEEeqnarray*}
where
\begin{gather}\nonumber
A =
\begin{bmatrix}
0.2229&	0.5637 \\
0.8708&	0.9984 
\end{bmatrix}, \, B_2  =
\begin{bmatrix}
0.5254 &	0.6644 \\
0.3872 &	0.9145 
\end{bmatrix}, \\
B_1 =\begin{bmatrix}
1 & 0  \\
0 &1  
\end{bmatrix},  \,  C=\begin{bmatrix}
1 & 0  \\
0 &1 \\
0&   0  \\
0&   0 
\end{bmatrix},\,
D=\begin{bmatrix}
0 & 0 \\0&0\\1 & 0 \\0&1
\end{bmatrix}. \nonumber
\end{gather}
\end{example}

\begin{figure}[t!]
\centering
\includegraphics[trim=0 0 0 0,width=1\columnwidth]{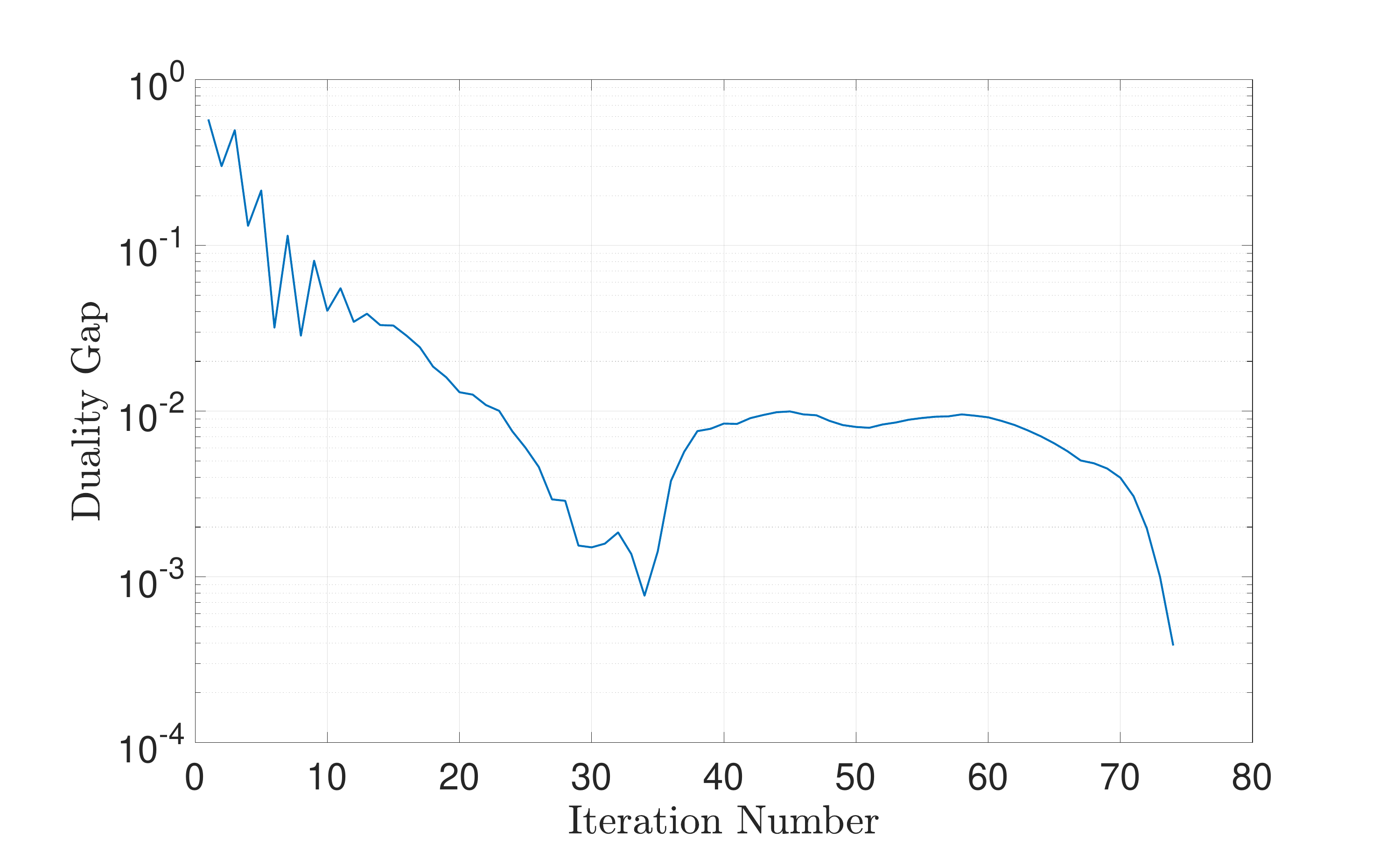}
\caption{Duality gap during iterations in Example 2}
\label{fig:cost2}
\end{figure}

Since all the parameters in $A$ and $B_2$ are uncertain with a variation of $\pm$20\%, a total of $2^{8}=256$ extreme systems need to be considered in the optimization.  
We performed 5 trials and recorded the total computation time with our approach, where the computation time is  5.3343 s, 5.6423 s, 5.1526 s, 5.8462 s, and 5.3149 s in these trials. In this case, the average computation time is 5.4581 s. Moreover, the change of the duality gap is shown in Fig.~\ref{fig:cost2}. At optimality, the following results are obtained, where
\begin{IEEEeqnarray*}{rCl}
W^* &=& \begin{bmatrix}
0.1608 &   0.0058 &   0.1673 &   0.0665\\
0.0058 &   0.0129 &   0.0327  &  0.0744\\
0.1673  &  0.0327  &  0.6892  &  0.2970\\
0.0665 &   0.0744 &   0.2970  &  1.2701
\end{bmatrix},\\
\mu^* &=& 0.0410,\\
K^*&=&	 \begin{bmatrix}
0.9643 & 2.1060 \\
0.2088 & 5.6843
\end{bmatrix},\\
\gamma^*&=&4.9411.
\end{IEEEeqnarray*}

\begin{figure}[t!]
\centering
\includegraphics[trim=0 0 0 0,width=1\columnwidth]{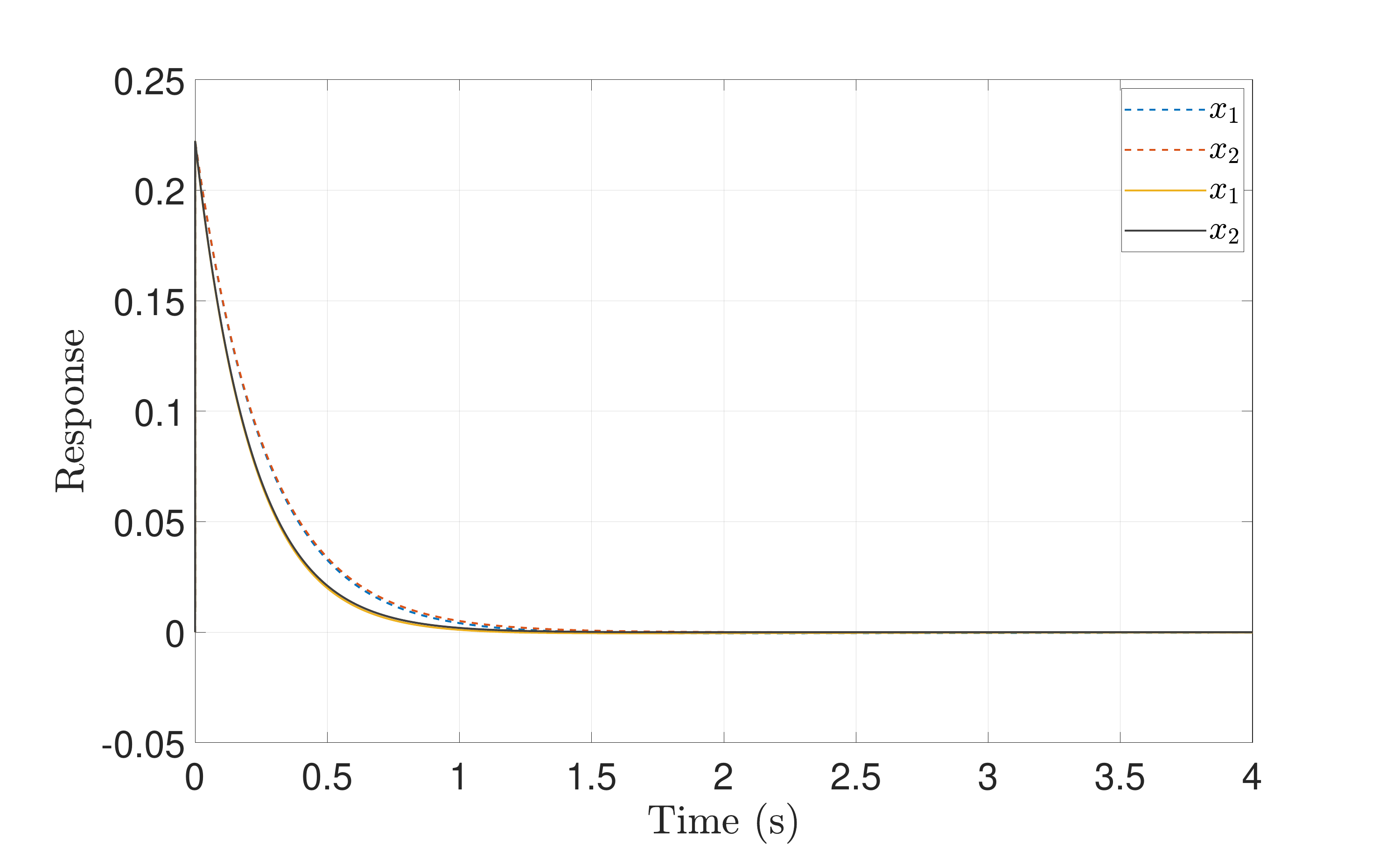}
\caption{System response in Example 2 (dashed line: extreme system; solid line: nominal system)}
\label{fig:response2}
\end{figure}

\begin{figure}[t!]
\centering
\includegraphics[trim=0 0 0 0,width=1\columnwidth]{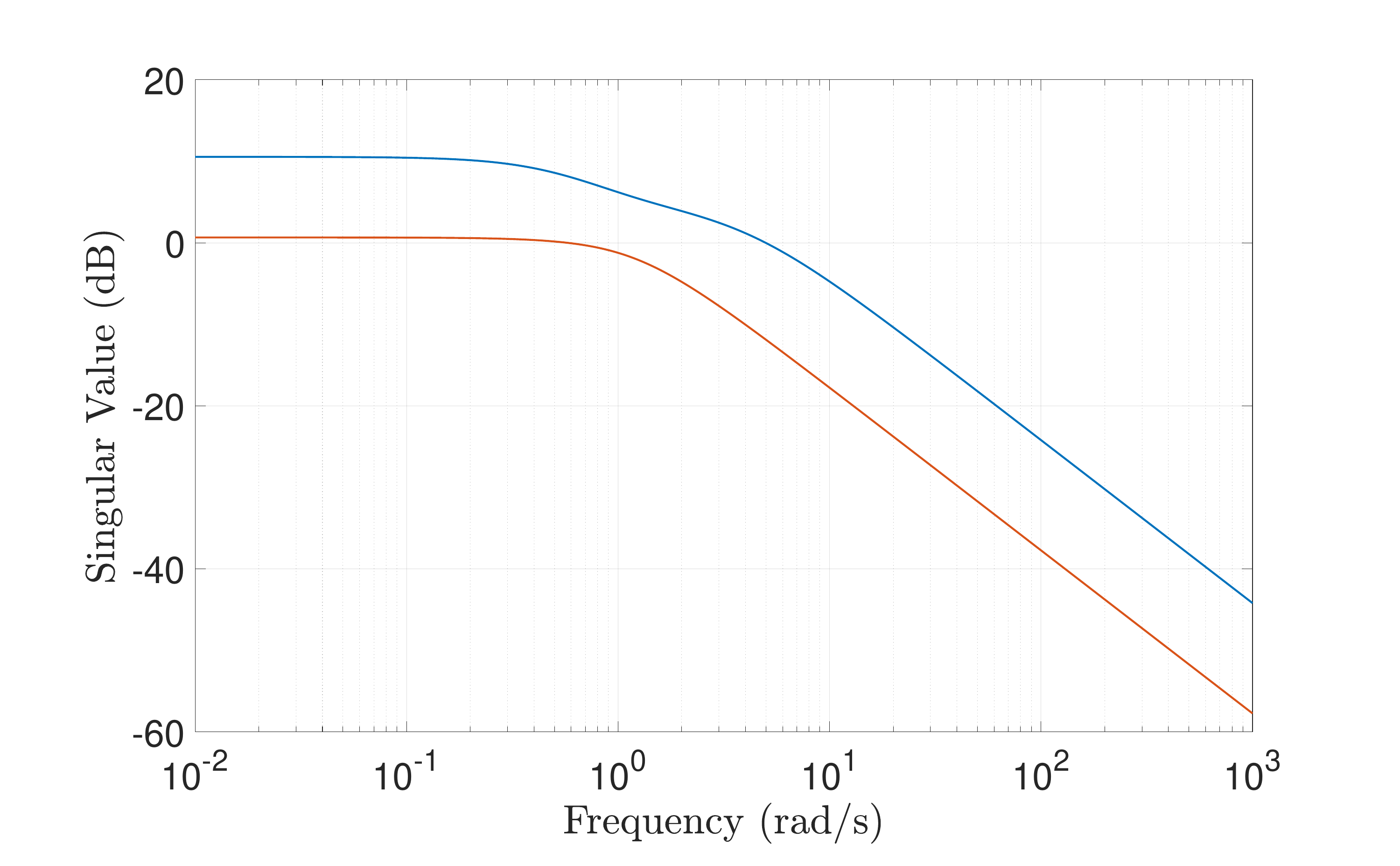}
\caption{Singular value diagram in Example 2}
\label{fig:sv2}
\end{figure}

In the simulation, $w$ is considered as a vector of the impulse disturbance. For illustration purposes, the simulation considers the nominal system and an extreme system with all the uncertain parameters reaching their lower bounds, then the responses of all state variables are shown in Fig.~\ref{fig:response2}. The dashed line shows the response of the extreme system, and the solid line indicates the response of the nominal system. It can be seen that for the extreme system, the closed-loop stability is suitably ensured despite the existence of parametric uncertainties. As clearly observed, the performance of the extreme system is slightly worse than the nominal system, but the difference of the dynamic response in terms of the extreme system and the nominal system is not significant. Hence, the robustness of the proposed approach is validated. Similarly, the singular value diagram of $H(s)$ is shown in Fig.~\ref{fig:sv2}, and the maximum singular value is given by 10.54 dB, which is equivalent to 3.3651 in magnitude, and it can be seen that it is bounded by $\gamma^*$.

Notice that the effectiveness of the proposed methodology in terms of computation can be more clearly demonstrated when there are a large number of extreme systems. Hence, a comparison study is carried out, where a well-established cutting-plane algorithm as presented in~\cite{peres1994optimal} is used. Note that this method has been demonstrated its effectiveness in solving a class of $\mathcal{H}_2$ and $\mathcal{H}_\infty$ problems in the parameter space. The numerical procedures are given below:
	
\noindent Step 1: Set $l=0$ and define the polytope $\mathscr{P}^0 \supseteq \mathscr{C}_U$.

\noindent Step 2: Solve the linear programming problem:
$(W^l,\mu^l)=\textup{argmax}\{\mu:(W, \mu)\in \mathscr{P}^l\}$.

\noindent Step 3: If $(W^l,\mu^l)\in \mathscr{C}_U$, $(W^l,\mu^l)$ is the optimal solution. Otherwise, generate a separating hyperplane and define $\mathscr{P}^{l+1}$. Set $l\leftarrow l+1$ and return to Step 2.

Essentially in the approach presented in~\cite{peres1994optimal}, a suitable polytope $\mathscr{P}^0$ is initialized such that $\mathscr{C}_U$ is a subset of $\mathscr{P}^0$. Then, the associated linear constraint is solved in the linear programming routine (which is conducted within  the polytope). For invoked nonlinear constraints, the cutting plane technique (also known as the outer linearization method) is adopted, where the satisfaction/violation condition of the nonlinear constraints is checked. If they are violated, half space (i.e., the cutting plane) will be generated and constructed for separation and update of the polytope in an iterative framework. Subsequently, the cutting planes are implemented as linear inequality constraints in the linear programming routine. In fact, this method makes an appropriate estimation of an unknown nonlinear set by iteratively involving a series of linear constraints without leading to infeasibility. However, with such a large number of extreme systems, the optimization process is unfortunately terminated with unsuccessful outcomes. The reason is that, in each iteration, a number of cutting planes could be generated and incorporated into the new linear programming routine, and these definitely lead to huge amount of computational efforts to solve this optimization problem. 

Additionally, some existing numerical solvers, such as Gurobi~\cite{gurobi} and SCS~\cite{SCS}, are not capable to solve this optimization problem.  
However, with an ADMM framework in our proposed development, the original optimization problem is decomposed into a series of manageable sub-problems that can be solved effectively. This is because that, in each iteration, the explicit solution to these sub-problems can be obtained very efficiently. In this case, our proposed approach alleviates the computational burden with the splitting scheme.

\begin{figure}[t!]
	\centering
	\includegraphics[trim=0 0 0 0,width=1\columnwidth]{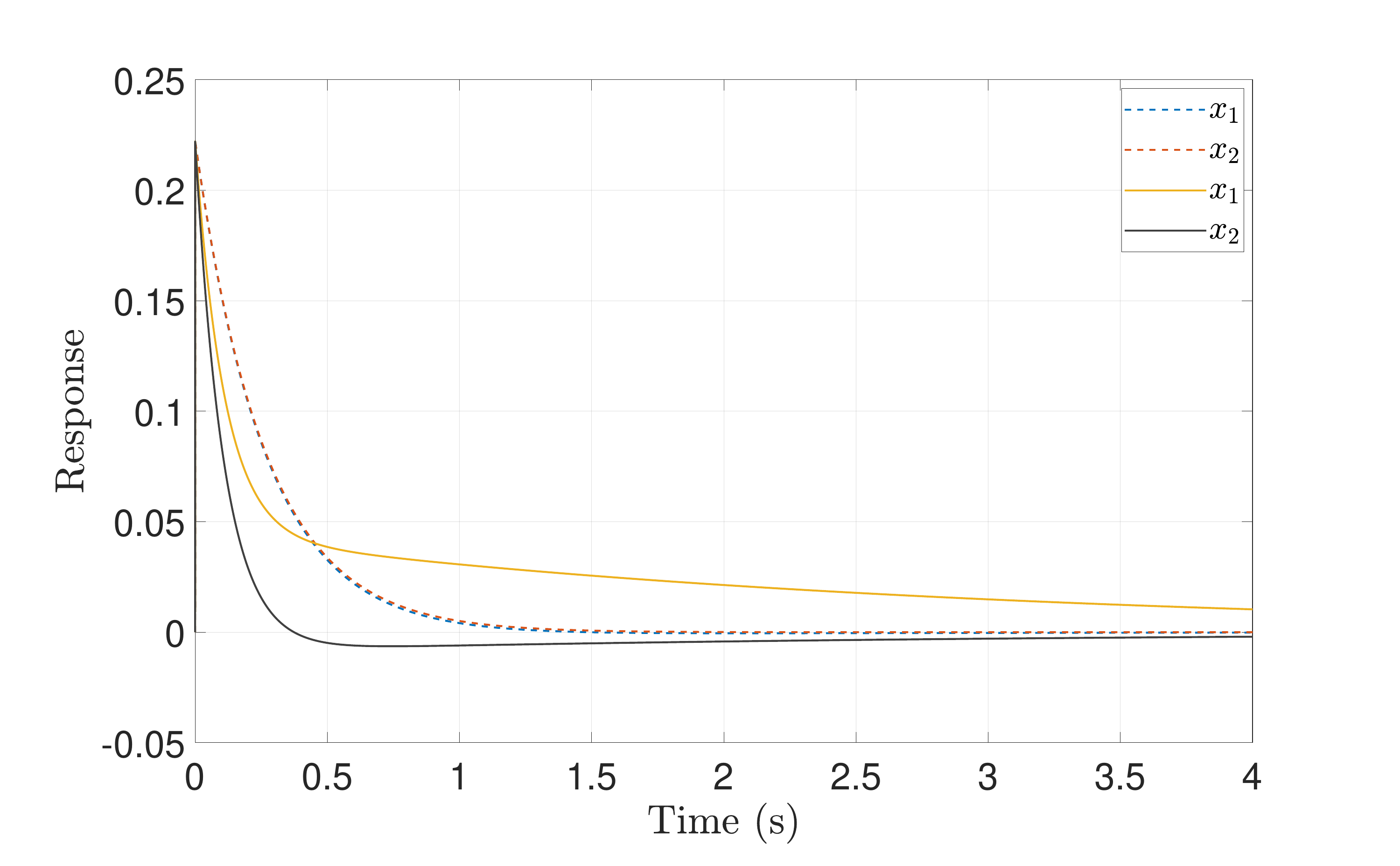}
	\caption{System response in Example 2 (dashed line: the proposed method; solid line: the LQR method)}
	\label{fig:response3}
\end{figure}

Furthermore, the LQR controller is used as a benchmark for comparison of the system performance. Note that for the LQR controller, it is determined by the algebraic Riccati equation. 
Then, the extreme system is used again in the comparison, and the system response is presented in Fig.~\ref{fig:response3}, where the dashed line represents the response with the proposed controller and the solid line indicates the response with the LQR controller. It is pertinent to note that for the purpose of a fair comparison, weighting matrices in the LQR are tuned such that the control inputs are at the same level as those attained in our proposed method.
Essentially our proposed method demonstrates superior performance over the LQR approach, which can be clearly observed from Fig.~\ref{fig:response3}. This is because parametric uncertainties are considered in our proposed method, and the disturbance attenuation is suppressed at its minimal level; however, the LQR does not take parametric uncertainties and disturbance attenuation into consideration.

%


\section{Conclusion}
In this work, the symmetric Gauss-Seidel ADMM algorithm is presented 
to solve the $\mathcal{H}_\infty$ guaranteed cost control problem, 
and the development and formulation of numerical procedures 
is given in detail
(with invoking a suitably interesting problem re-formulation based on the Schur complement).
Through a parameterization technique
(where the stabilizing controllers 
are characterized 
by an appropriate convex parameterization 
which is described and established analytically
in our work here), 
the robust stability and performance can be suitably achieved 
in the presence of parametric uncertainties. 
An upper bound of all feasible $\mathcal{H}_\infty$ performances 
is minimized over the uncertain domain, 
and the minimum disturbance attenuation level 
is obtained through the optimization. 
Furthermore, the algorithm is evaluated based on two suitably appropriate illustrative examples, 
and the simulation results successfully reveal 
the practical appeal of the proposed methodology in terms of computation, 
and also clearly validate the results on robust stability and performance.
This work can be further extended to the mixed $\mathcal{H}_2$/$\mathcal{H}_\infty$ control problem, which aims to balance the trade-off between performance and robustness. Particularly, given an prescribed $\mathcal{H}_\infty$ attenuation level, the objective is to seek the optimal $\mathcal{H}_2$ control gain, where the nominal performance index is minimized with the imposed pertinent constraints.

\bibliographystyle{IEEEtran}
\bibliography{IEEEabrv,Reference}

\end{document}